\documentclass[a4paper,11pt]{article}
\usepackage{amsmath,amssymb}

\newcommand{\dd}{{\rm d}} \newcommand{\e}{{\rm e}} \newcommand{\im}{{\rm i}}
\newcommand{\Tr}{{\rm Tr}} \renewcommand{\Im}{{\rm Im}} \renewcommand{\Re}{{\rm Re}}
\newcommand{\bbbone}{{\mathchoice {\rm 1\mskip-4mu l} {\rm 1\mskip-4mu l}
{\rm 1\mskip-4.5mu l} {\rm 1\mskip-5mu l}}}
\newcommand{\esssup}{\mathop{\rm ess\ sup}} \renewcommand{\rho}{\varrho}
 
 \newcommand{\Dscr}{{\mathcal{D}}}
\newcommand{\Escr}{{\mathcal E}} \newcommand{\Fscr}{{\mathcal F}}
 \newcommand{\Hscr}{{\mathcal{H}}}
\newcommand{\Iscr}{{\mathcal{I}}} \newcommand{\Lscr}{{\mathcal{L}}}
 \newcommand{\Mscr}{{\mathcal{M}}}

\newcommand{\Rscr}{{\mathcal{R}}} \newcommand{\Sscr}{{\mathcal{S}}}
\newcommand{\Tscr}{{\mathcal{T}}} \newcommand{\Uscr}{{\mathcal{U}}}
 \newcommand{\Wscr}{{\mathcal{W}}}
\newcommand{\Xscr}{{\mathcal{X}}} 
\newcommand{\Zscr}{{\mathcal{Z}}} \newtheorem{theorem}{Theorem}[section]
\newtheorem{assumption}{Assumption}[section]

\newtheorem{definition}{Definition}[section]

\begin{document}

\title{On a class of stochastic differential equations \\ used in quantum optics
\footnote{To appear in {\sl Rendiconti del Seminario Matematico e Fisico di Milano}}}
\author{A.~Barchielli\dag\ and  F.~Zucca\ddag\\
{\small \dag\ {\sl Dipartimento di Matematica, Politecnico di Milano,}}\\
{\small{\sl Piazza Leonardo da Vinci 32, I-20133 Milano, Italy}} \\
{\small and  {\sl Istituto Nazionale di Fisica Nucleare, Sezione di Milano}} \\
{\small \ddag\ {\sl Dipartimento di Matematica, Universit\`a degli Studi di Milano, }}\\
{\small {\sl Via Saldini 50, I-20133 Milano, Italy}}}
\date{}
\maketitle

\numberwithin{equation}{section}

\begin{abstract}
Stochastic differential equations for processes with values in Hilbert spaces
are now largely used in the quantum theory of open systems. In this work we
present a class of such equations and discuss their main properties; moreover,
we explain how they are derived from purely quantum mechanical models, where the
dynamics is represented by a unitary evolution in a Hilbert space, and how they
are related to the theory of continual measurements. An essential tool is an
isomorphism between the bosonic Fock space and the Wiener space, which allows to
connect certain quantum objects with probabilistic ones.
\end{abstract}

\section{Classical stochastic differential equations}

Recently there was an increasing use of stochastic differential equations (SDEs) in
quantum optics and more generally in the theory of quantum open systems
\cite{1G}--\cite{6B}. In particular SDEs are used for theoretical purposes \cite{2B,3B,6B}
in the theory of quantum continuous measurements \cite{6B,BLP}, when a system is
continually monitored in time, and for numerical simulations of quantum master equations
\cite{4G}, which are mathematical models for the dynamics of quantum open systems. In this
work we want to describe a restricted class of such equations, to discuss their main
properties and to show how they can be derived from fully quantum mechanical models
\cite{7B,BPag1}.

We consider only SDEs of diffusive type, when the noises are Wiener processes; for the
theory of It\^o's stochastic integrals in Hilbert spaces we refer to \cite{DPZ}. Let
$\big(\Omega, (\Fscr_t), \Fscr, P\big)$ be a stochastic basis satisfying the usual
hypotheses; in particular, $(\Omega, \Fscr,P)$ is a probability space, which we assume to
be standard, $(\Fscr_t)_{t\geq 0}$ is a filtration, i.e. $\Fscr_t$ is a
sub--$\sigma$--algebra of $\Fscr$ and $\Fscr_s
\subset \Fscr_t$ for $s \leq t$; we also assume $\Fscr = \bigvee_{t \geq 0} \Fscr_t$. Let
$W_k(\cdot)$, $k=0,1,\ldots$, be a sequence of continuous versions of adapted, standard,
independent, Wiener processes with increments independent from the past.

Let $\Hscr$ be a separable complex Hilbert space and let us consider the
following linear SDE with ``multiplicative noise'' for an $\Hscr$--valued
process $\psi_t$:
\begin{equation} \label{sd1}
\left\{\begin{array}{rcl}
\dd \psi_t & = &  {\displaystyle\sum_{j=0}^\infty} R_j(t) \psi_{t}\,
\dd W_j (t) - \im K(t) \psi_{t}\, \dd t \\
\psi_0 & = & \psi
\end{array}\right.
\end{equation}
where the coefficients $R_j(t)$, $K(t)$ are operators in $\Hscr$ satisfying some
assumptions to be discussed below and the initial condition $\psi$ is an
$\Hscr$--valued $\Fscr_0$--measurable random variable. Let us state now some
different sets of assumptions about the coefficients.

\begin{assumption} \label{A1}
The coefficients are time dependent, bounded operators on $\Hscr$ satisfying:
\item{\rm (i)} $\sum_{j=0}^\infty R^*_j(t) R_j(t)$ is strongly convergent in
$\Lscr (\Hscr)$ for all $t$;
\item{\rm (ii)} the functions $t \mapsto R(t)$, $t
\mapsto K(t)$ are strongly measurable, where $R(t)$ is the operator from $\Hscr$
to $\Hscr \otimes l^2$ defined by $ \left \langle x \otimes c | R(t) y \right
\rangle = \sum_{j=0}^\infty \overline{c_j} \left \langle x | R_j(t) y
\right \rangle$, $\forall x,y \in \Hscr$, $\forall c \in l^2$;
\item{\rm (iii)}
$\displaystyle
\forall T \in \mathbb R_+, \qquad \esssup_{t \leq T} \bigg(
\bigg\|\sum_{j=0}^\infty R^*_j(t) R_j(t)\bigg\|^2 + \|K(t)\| \bigg ) < \infty\,.
$
\end{assumption}

\begin{assumption} \label{A2}
The coefficients are time independent, (in general) unbounded operators, defined
on a dense domain $\Dscr \subset \Hscr$; we assume that
\[
\sum_{j=0}^\infty \|R^*_j x\|^2_{{}_\Hscr} < \infty, \qquad \forall x \in \Dscr^*,
\]
where $\Dscr^*$ is a core for the adjoint operators $K^*$ and $R_j^*$. Moreover
we assume the following dissipativity condition to hold:
$$
\sum_{j=0}^\infty \|R_j x\|^2_{{}_\Hscr} + 2 \Im \left \langle x | K x
\right \rangle_{{}_\Hscr}
\leq c \|x\|^2_{{}_\Hscr}\,, \qquad \forall x \in \Dscr.
$$
\end{assumption}

\begin{assumption} \label{A3}
The coefficients are time independent, (in general) unbounded operators, defined
on a dense domain $\Dscr \subset \Hscr$, satisfying
$$
\sum_{j=0}^\infty \|R_j x\|^2_{{}_\Hscr} < \infty\,, \qquad \forall x \in \Dscr;
$$
moreover, $-\im K$ is the generator of an analytic operator semigroup $T(t)$.
\end{assumption}

Let us remark that, in the time independent case Assumption \ref{A1} is stronger
than Assumption \ref{A2} and indeed Theorem \ref{eu1} will show that, under
Assumption \ref{A1}, we can obtain stronger results than under Assumption
\ref{A2}. We note also that, for a linear equation, the so called ``Lipschitz
conditions'' \cite{DPZ} are equivalent to the boundness of the coefficients.

\begin{definition}\label{solution}{\rm
We consider three different kinds of solution:
\vskip 8 pt
\noindent {\rm (i)} a (topological) \emph{strong solution} of eq.~\eqref{sd1} is
an $\Hscr$--valued process such that $\forall t
\in \mathbb R_+$ the following equation holds almost surely (a.s.)
\begin{equation}\label{sde1}
\psi_t  =  \psi+\sum_{j=0}^\infty \int_0^t R_j(s) \psi_{s} \,\dd W_j (s)
- \im \int_0^t K(s) \psi_{s}\, \dd s\,;
\end{equation}
{\rm (ii)} a (topological) \emph{weak solution} of eq.~\eqref{sd1} is an $\Hscr$--valued
process such that $\forall t\in
\mathbb R_+$ and $\forall \chi \in \Dscr^*$ the following equation holds a.s.
\begin{equation}\label{sde2}
\left \langle \chi | \psi_t \right \rangle_{{}_\Hscr}  =
\left \langle \chi | \psi \right \rangle_{{}_\Hscr} +
\sum_{j=0}^\infty \int_0^t \left \langle \left. R^*_j(s) \chi \,\right|
\psi_{s} \right \rangle_{{}_\Hscr}\, \dd W_j (s) - \im
\int_0^t \left \langle K^*(s) \chi | \psi_{s} \right \rangle_{{}_\Hscr}\, \dd s\,;
\end{equation}
{\rm (iii)} a (topological) \emph{mild solution} of eq.~\eqref{sd1} with Assumption
\ref{A3} is an $\Hscr$-valued process such that $\forall t
\in \mathbb R_+$ the following equation holds a.s.
\begin{equation}\label{sde3}
\psi_t  =  T(t)\psi+ \sum_{j=0}^\infty \int_0^t T(t-s)R_j \psi_{s}\, \dd W_j (s)\,,
\end{equation}
where $T(s)$ is the analytic semigroup operator generated by $- \im K$.}
\end{definition}

For more details about these definitions see \cite{DPZ}. It is clear that a
strong solution is also a weak solution and it is possible to show that a weak
solution is also a mild solution (see \cite{DPZ}, Propositions 6.2, 6.3, 6.4 and
Theorem 6.5). In the following theorem we state some results about the
existence, the uniqueness and the Markov properties of the solution of
eq.~\eqref{sd1} under different hypotheses.

\begin{theorem} \label{eu1}
{\rm (i)} Under Assumption \ref{A1}, eq.~\eqref{sd1} admits a unique strong
solution;
{\rm (ii)} under Assumption \ref{A2}, eq.~\eqref{sd1} admits a unique
weak solution;
{\rm (iii)} under Assumption \ref{A3}, eq.~\eqref{sd1} admits a
mild solution and this is also a weak solution. All these solutions are strong
Markov processes.
\end{theorem}
\noindent {\it Proof.} The statement (i) is Theorem 1.1 of \cite{BPZ}; for the
proof of (ii) see Theorem 1 of \cite{Hol1} and for the proof of (iii) see
Theorems 6.19 and 6.22 of \cite{DPZ}. For the Markov property see \cite{Hol1}
Theorem 1, \cite{DPZ} Theorems 9.8 and 9.14. \hfill{$\square$}

\smallskip

We are not interested in eq.~\eqref{sd1} in general, but only when $\|
\psi_t\|^2_{{}_\Hscr} $ is a martingale and can be interpreted as a probability
density with respect to $P$. Indeed, if we have a positive martingale
$\|\psi_t\|^2_{{}_\Hscr}$ with $\mathbb E_{{}_P} \left [\|\psi_t\|^2_{{}_\Hscr}
\right ]=1$ ($\mathbb E_{{}_P}$ denotes the expectation with respect to $P$),
then
\begin{equation} \label{sd2}
\widehat P (A) := \mathbb E _{{}_P}\left[ 1_A \|\psi_t\|^2_{{}_\Hscr}\right],
\qquad \forall A \in \Fscr_t\,, \ \forall t \in \mathbb R_+,
\end{equation}
defines a new probability law $\widehat P$ on $(\Omega, \Fscr)$. To obtain this
we need a further assumption.

\begin{assumption}\label{A4}
For all $t \in \mathbb R_+$, the coefficients satisfy the following operator equality
$$
-\im \left (K^*(t) - K(t) \right ) =\sum_{j=0}^\infty R^*_j(t)R_j(t),
$$
(where the sum is convergent in the strong operator topology) and the initial
condition is a square integrable, normalized random variable, i.e.\ $\mathbb
E_{{}_P} \left [\|\psi\|^2_{{}_\Hscr} \right ]=1$.
\end{assumption}

\begin{theorem}\label{msc}
Under Assumptions \ref{A1} and \ref{A4} the square norm
$\|\psi_t\|^2_{{}_\Hscr}$of the solution is a (positive) martingale and $\mathbb
E_{{}_P} \left [\|\psi_t\|^2_{{}_\Hscr} \right ]=1$. Moreover, under the law
$\widehat P$ defined by eq.~\eqref{sd2}, the processes
\begin{equation}\label{sd3}
\widehat W_k(t):= W_k(t)-2 \int_0^t \Re \, \widehat m_k(s) \dd s, \qquad k=0,1,\ldots,
\end{equation}
where
\begin{equation}\label{sde4}
\widehat m_k(t) :=
\frac{ \left \langle \psi_{t} | R_k(t) \psi_{t} \right
\rangle_{{}_\Hscr}}{\|\psi_{t}\|^2_{{}_\Hscr}}\,,
\end{equation}
are independent, adapted, standard Wiener processes.
\end{theorem}

For the proof of this theorem see \cite{BPZ}, Theorem 1.2 and Proposition 1.1.
The unbounded case is more difficult and it is treated in \cite{Hol1} under an
additional condition called of  hyperdissipativity.

Let us introduce now the process
\begin{equation}\label{sde5}
\widehat \psi_t :=
\frac{1}{\|\psi_t\|_{{}_\Hscr}} \,\exp \bigg ( - \im \sum_{k=0}^\infty \int_0^t \left
( \Im \widehat m_k(s)\right) \left(\dd W_k(s) - \Re\, \widehat m_k(s)\,
\dd s\right) \bigg )\,\psi_t\,.
\end{equation}
Note that $\|\widehat \psi_t\|_{{}_\Hscr}^2=1$; indeed, $\exp(\cdots)$ is  a
stochastic phase with no particular meaning, which we introduce only in order to
simplify the form of the stochastic differential of $\widehat \psi_t$. Note also
that $\widehat m_k(t) = \left \langle \widehat \psi_t \left| R_k(t) \widehat
\psi_t \right.\right \rangle$.

\begin{theorem}\label{eu2}
Under Assumptions \ref{A1} and \ref{A4}, in the stochastic basis
$(\Omega,(\Fscr_t), \Fscr, \widehat P)$, $\widehat \psi_t$ satisfies the non
linear SDE
\begin{equation} \label{sd4}
\dd \widehat \psi_t = - \im \widehat K\left(t, \widehat \psi_t\right)\dd t +
\sum_{k=0}^\infty \widehat R_k \left(t, \widehat \psi_t\right) \dd \widehat W_k(t),
\end{equation}
where, $\forall f \in \Hscr$,
\begin{multline} \label{sd5}
\widehat K(t, f) := \frac{1}{2}(K(t)+K^*(t))f - \frac{\im}{2}
\sum_{k=0}^\infty
\left ( \frac{\left \langle f | R_k(t) f \right
\rangle_{{}_\Hscr}}{\|f\|^2_{{}_\Hscr}}\, R^*_k(t)
- \frac{\left \langle f | R^*_k(t) f \right
\rangle_{{}_\Hscr}}{\|f\|^2_{{}_\Hscr}} \,R_k(t)
\right )f + {} \\
{}-\frac{\im}{2} \sum_{k=0}^\infty
\left ( R^*_k(t) - \frac{\left \langle f | R^*_k(t) f \right
\rangle_{{}_\Hscr}} {\|f\|^2_{{}_\Hscr}} \right ) \left ( R_k(t)
- \frac{\left \langle f | R_k(t) f \right \rangle_{{}_\Hscr}}
{\|f\|^2_{{}_\Hscr}} \right ) f\,,
\end{multline}
and
\begin{equation} \label{sd6}
\widehat R_k (t, f) := \left ( R_k(t)
- \frac{\left \langle f | L_k(t) f \right \rangle_{{}_\Hscr}}
{\|f\|^2_{{}_\Hscr}} \right ) f\,.
\end{equation}
Moreover $\widehat \psi_t$ is the unique strong solution of eq.~\eqref{sd4} and
it is a strong Markov process.
\end{theorem}
\noindent {\it Proof.} The first part of this theorem is essentially Theorem 1.3
of \cite{BPZ}. Moreover it is easy to show that the coefficients $\widehat K$
and $\widehat R_k$ satisfy the local Lipschitz conditions (see \cite{DPZ},
p.~198); hence $\widehat \psi_t$ is the unique strong solution of
eq.~\eqref{sd4}. Again the Markov property follows from \cite{DPZ} Theorems 9.8
and 9.14. \hfill{$\square$}

\smallskip

Let us stress that the Markov property allows to associate to the process
$\widehat \psi_t$ a transition function satisfying a Chapman--Kolmogorov
equation (see \cite{DPZ} Section 9.2.1 and, in particular, Corollary 9.9);
similar considerations apply to the solution of the linear SDE.

Examples of equations of the types \eqref{sd1} and \eqref{sd4} were introduced
in the theory of quantum continual measurement in Ref.~\cite{2B}.

In the statistical formulation of quantum mechanics the states of a quantum
system are represented by \emph{statistical operators}. Let $\Tscr(\Hscr)$ be
the \emph{trace class} on $\Hscr$ and $\Sscr(\Hscr)$ denote the set of
statistical operators, i.e.
\begin{align*}
\Tscr(\Hscr) &:= \left\{ \rho \in \Lscr(\Hscr) : \Tr_{{}_\Hscr}
\left\{ \sqrt{\rho^*\rho} \right\} <\infty \right\},
\\
\Sscr(\Hscr) &:= \left\{ \rho \in \Tscr(\Hscr) : \rho^*=\rho,\,
\rho\geq 0,\, \Tr_{{}_\Hscr} \left\{\rho\right\}=1 \right\}.
\end{align*}
Equipped with the norm $\|\rho\|_1 :=  \Tr_{{}_\Hscr}\left\{
\sqrt{\rho^*\rho}\, \right\}$, $\Tscr(\Hscr)$ becomes a Banach space whose dual
is $\Lscr(\Hscr)$.

\begin{theorem}\label{ME}
Let $\psi_t$ be the solution of eq.~\eqref{sd1} with Assumptions \ref{A1} and
\ref{A4}. Then, the formula
\begin{equation}\label{sd8}
\Tr_{{}_\Hscr}\left\{\rho_t\, a\right\} = \mathbb E_{{}_P} \left [
\left \langle \psi_t | a \psi_t
\right \rangle_{{}_\Hscr} \right ], \qquad \forall a \in \Lscr(\Hscr),
\end{equation}
defines a family $\{\rho_t,\ t\geq 0\}$ of statistical operators on $\Hscr$ with
$\rho_0=\rho$, where $
\Tr_{{}_\Hscr}\left\{\rho\, a\right\} := \mathbb E_{{}_P}\left [
\left \langle \psi | a \psi
\right \rangle_{{}_\Hscr} \right ]$, $\forall a \in \Lscr(\Hscr)$;
eq.~\eqref{sd8} is equivalently written as
\begin{equation}\label{sd12}
\Tr_{{}_\Hscr}\left\{\rho_t\, a\right\} = \mathbb E_{{}_{\widehat P}}
\left [ \left \langle \widehat \psi_t \left| a \widehat \psi_t \right.
\right \rangle_{{}_\Hscr} \right ], \qquad \forall a \in \Lscr(\Hscr).
\end{equation}
Moreover, $\rho_t$ is the unique solution of the integral equation
\begin{equation}\label{sd9}
\rho_t = \rho + \int_0^t \Lscr(s)_{*} [\rho_s]\, \dd s,
\end{equation}
where $\Lscr(s)_{*}$ is the preadjoint operator of $\Lscr(s) \in \Lscr (\Lscr(\Hscr
))$ defined by
\begin{equation}\label{sd11}
\Lscr(t)[a] := - \frac{\im}{2} \,[K(t)+K^*(t),a] +\frac{1}{2}
\sum_{i=0}^\infty \left( R_i^*(t)\big[a,R_i(t)\big]+
\big[ R^*_i(t),a\big]R_i(t)\right), \qquad \forall a\in \Lscr(\Hscr),
\end{equation}
where $[a,b]:=ab-ba$.
\end{theorem}

For the proof of this theorem see \cite{BPZ} Theorem 1.2 and Section 3.

In the theory of open quantum systems eq.~\eqref{sd9} is usually written as $\dd
\rho_t/ \dd t = \Lscr(t)_*[\rho_t]$ (if $\Lscr(t)$ is sufficiently smooth in
time), it is called a \emph{master equation} and  describes the evolution of the
statistical operator representing the state of the system.

Formula \eqref{sd12} connects the solution $\widehat \psi_t$ of the SDE
\eqref{sd4} with the solution $\rho_t$ of the master equation \eqref{sd9}; this
fact is at the basis of stochastic methods for numerical simulations of
solutions of master equations \cite{4G}.

For what concerns the theory of continual measurements, we shall see in Section
3 that the processes $W_k(t)$, or functionals of them, represent the output of
the measurement and that the physical law of this output is given by the
probability $\widehat P$ defined in eq.~\eqref{sd2}. In particular, formulae for
all the moments (under $\widehat P$) of the processes $W_k(t)$ can be given. For
instance, from Theorems \ref{msc} and \ref{ME} we get
$
\mathbb E_{{}_{\widehat P}} \left[ W_k(t)\right] = 2 \int_0^t \Re\,
\mathbb E_{{}_{\widehat P}} \left[ \widehat m_k(s)\right] \dd s$ and $
\mathbb E_{{}_{\widehat P}} \left[ \widehat m_k(s)\right]= \Tr_{{}_\Hscr}
\left\{ \rho_s \, R_k(s) \right\}
$; therefore, the mean value of the output $W_k(t)$ is connected to the solution
of the master equation by
\begin{equation}\label{sdx}
\mathbb E_{{}_{\widehat P}} \left[ W_k(t)\right] =  \int_0^t
\Tr_{{}_\Hscr} \left\{ \rho_s \bigl( R_k(s) + R^*_k(s) \bigr) \right\} \dd s.
\end{equation}
The computation of the moments of the second order is a little bit more
involved; one needs to go through the characteristic functional \cite{BLP,QO} or
to exploit the Markov properties of the solution of eq.~\eqref{sd4}. Let
$\Uscr(t,s)$ be the evolution operator (from time $s$ to $t$) associated with
eq.~\eqref{sd9}, i.e.\ $\rho_t = \Uscr(t,0)[\rho]$; then, the result is (cf
\cite{QO} eq.~(5.20))
\begin{multline}\label{200}
\mathbb E_{{}_{\widehat P}} \left[ W_i(t_1) W_j(t_2) \right] = t_1 \wedge t_2 +{} \\
{}+ \int_0^{t_1} \dd s_1 \int_0^{t_2\wedge s_1} \dd s_2 \, \Tr_{{}_\Hscr}
\left\{ \Bigl( R_i(s_1) + R_i^*(s_1) \Bigr) \Uscr(s_1,s_2)
\left[ R_j(s_2) \rho_{s_2} + \rho_{s_2} R_j^*(s_2) \right] \right\} +{} \\
{}+ \int_0^{t_2} \dd s_2 \int_0^{t_1\wedge s_2} \dd s_1 \, \Tr_{{}_\Hscr}
\left\{ \left( R_j(s_2) + R_j^*(s_2) \right) \Uscr(s_2,s_1) \Bigl[ R_i(s_1)
\rho_{s_1} + \rho_{s_1} R_i^*(s_1) \Bigr] \right\} .
\end{multline}

Also the random vectors $\widehat \psi_t$ have a physical interpretation
\cite{2B,3B,5C,6B}. Indeed, it can be shown that the random orthogonal
projection $\left| \widehat \psi_t \right\rangle \left \langle \widehat \psi_t
\right|$ represents the state at time $t$ one has to attribute to the quantum
system if the trajectory of all the processes $W_j$ is known up to time $t$; in
this sense $\left| \widehat \psi_t \right\rangle \left \langle \widehat \psi_t
\right|$ represents the \emph{a posteriori} state of the system at time $t$.
Then, $\rho_t = \mathbb E_{{}_{\widehat P}}\left[ \left| \widehat \psi_t
\right\rangle \left \langle \widehat \psi_t \right|\right]$ is the state of the
system when no information on the outputs is taken into account; therefore,
$\rho_t$ is called the \emph{a priori} state at time $t$.

Note that no physical quantity, such as $\rho_t$, $\left| \widehat \psi_t
\right\rangle \left \langle \widehat \psi_t \right|$ and $\widehat P$, depends
on the arbitrary phase introduced in eq.~\eqref{sde5}.

\section{The Fock and Wiener spaces}

In order to be able to show how the stochastic equations of the previous section
are related to the quantum dynamics, we need to recall some facts about Fock and
Wiener spaces.

Let $\Zscr$ be a separable complex Hilbert space and let us denote by $\Gamma$
the symmetric Fock space over $\Xscr:=\Zscr\otimes L^2({\mathbb R}) \simeq
L^2({\mathbb R};\Zscr)$, which is defined by
\begin{equation} \label{d.1}
\Gamma := {\mathbb C} \oplus {\sum_{n=1}^\infty}^{{}_{{}_{{}_{\!\oplus}}}}
\left( \Xscr^{\otimes n}\right)_{\rm s}\, ,
\end{equation}
where $\left( \Xscr^{\otimes n}\right)_{\rm s}$ denotes the symmetric part of
the tensor product $\Xscr^{\otimes n}$ and is called the $n$-particle space
(\cite{Par} Sects.~17 and 19). The vectors $ E(f) := \left( 1, f, \frac{1
}{\sqrt{2!}}\, f\otimes f, \ldots,  \frac{1 }{\sqrt{n!}}\, f^{\otimes n}, \ldots
\right)$,  $f\in \Xscr$, are called \emph{exponential vectors} (\cite{Par}
p.~124); let us denote by $\Escr$ the linear span of the exponential vectors.

\begin{theorem}\label{t1}
The set of all the exponential vectors is linearly independent and total in
$\Gamma$; moreover, for $f,\,g\in \Xscr$ we have $
\langle E(g) | E(f) \rangle_{{}_\Gamma} = \exp \langle g|f\rangle_{{}_\Xscr}$.
Let us fix $g\in \Xscr$ and $V \in \Uscr(\Xscr)$ (unitary operators on $\Xscr$);
there exists a unique unitary operator $\Wscr(g;V)$ \emph{(Weyl operator)} such
that
\begin{equation}\label{d.4}
\Wscr(g;V) E(f):= \exp\left\{ -\frac{1 }{2} \|g\|^2_{{}_\Xscr} - \langle
g|Vf\rangle_{{}_\Xscr} \right\} E(Vf+g), \qquad \forall f\in \Xscr.
\end{equation}
These operators satisfy the multiplication rules $\big( g_i \in \Xscr$, $V_i
\in \Uscr(\Xscr)\big)$
\begin{equation}\label{d.5}
\Wscr(g_2;V_2) \Wscr(g_1;V_1) =\exp \left [\im\, \Im \langle V_2g_1|g_2
\rangle_{{}_\Xscr} \right] \Wscr( V_2g_1+g_2;\,V_2V_1)\,.
\end{equation}\end{theorem}

Total means that $\Escr$ is dense in $\Gamma$; for the proof see \cite{Par}
p.~124,  Proposition 19.4 p.~126, pp.~134--135. Note that
$\|E(f)\|_{{}_\Gamma}^2 = \exp \|f\|_{{}_\Xscr}^2$; by normalizing the
exponential vectors one obtains the \emph{coherent vectors} $e(f) := \exp\left[
- \frac{1}{2} \|f\|^2_{{}_\Xscr} \right]\, E(f)$. The vector $e(0)\equiv E(0)
\equiv (1,0,0,\ldots)$ is called the \emph{vacuum}.

In the following we shall need particular Weyl operators. Let $V\in
\Uscr(\Xscr)$ be such that $\big(Vf\big)(t) = V(t) f(t)$, $\forall f\in \Xscr
\simeq L^2( \mathbb R;\Zscr)$; then $V(t)\in \Uscr(\Zscr)$. Let $g\in L^2_{\rm
loc}( \mathbb R;\Zscr)$, i.e.\ $\int_s^t \|g(\tau)\|^2_{{}_\Zscr} \,\dd \tau <
+\infty$, $\forall s,t\in \mathbb R$, $s<t$. Then, we set
\begin{equation}\label{d.2}
\Wscr_t(g;V) := \Wscr\left(\chi_{[0,t]}g;\,\chi_{[0,t]}(V-\bbbone)+ \bbbone \right).
\end{equation}

The Fock space has a \emph{continuous tensor product structure}, i.e.\ $\forall
s,t \in \mathbb R$, $s\leq t$, $\Gamma = \Gamma^s_{-\infty} \otimes \Gamma^t_s
\otimes \Gamma_t^{+\infty}$, where $\Gamma^b_a$ is  the symmetric Fock space
over $L^2\big((a,b);\Zscr\big)$, $-\infty \leq a < b\leq +\infty$.
Correspondingly one has $E(f) = E\left(\chi_{(-\infty,s)}f \right) \otimes
E\left(\chi_{(s,t)}f \right) \otimes E\left(\chi_{(t,+\infty)}f \right)$; we
also denote by $\Escr_a^b \subset \Gamma_a^b$ the linear span of the exponential
vectors of the type $E\left(\chi_{(a,b)} \cdot\right)$ (\cite{Par} Proposition
19.6 p.~127 and pp.~179--180). Note that $\Wscr_t(g;V)$ leaves invariant
$\Gamma_0^t$. By setting $S(t)E(f):= E(f_t)$, $f_t(s):=f(s-t)$, one defines a
strongly continuous unitary group of shift operators on $\Gamma$; moreover, one
has $S(t) \Gamma_a^b = \Gamma_{a+t}^{b+t}$.

Let $\{ e_j,\, j=0,1,\ldots\}$ be a c.o.n.s.\ in $\Zscr$, $g\in \Xscr$ and let
us define on $\Escr$ the operators
\begin{eqnarray}\label{d.6}
a(g)E(f) := \langle g|f\rangle_{{}_\Xscr}\, E(f)\,, \quad \forall
f\in \Xscr, \qquad &{}& A_j(t):= a\left(\chi_{[0,t]}(\cdot)e_j\right),\\
a^\dagger(g) E(f) := \frac{\dd \ }{\dd \varepsilon} \,
E\left( f+\varepsilon g \right)\Big|_{\varepsilon=0}\,, \quad
\forall f\in \Xscr, \qquad &{}& A_j^\dagger(t):=
a^\dagger\left(\chi_{[0,t]}(\cdot)e_j\right);\label{d.7}
\end{eqnarray}
let us note that one has $\langle a^\dagger(g) E(f_1)| E(f_2)\rangle_{{}_\Gamma}
= \langle E(f_1)| a(g) E(f_2) \rangle_{{}_\Gamma}$. The operator $a^\dagger(g)$
turns out to be linear in its argument and $a(g)$ antilinear. The domains of
these operators can be suitably extended; there exists a core $\Dscr$ for all
these operators where all their products are well defined and such that $\Dscr
\supset \Escr$, $\Dscr \supset (\Xscr^{\otimes n})_{\rm s}$ for every $n$. Then,
$A_j(t)$ and $a(g)$ send the $n$-particle space into the $(n-1)$--one and
$A_j^\dagger(t)$ and $a^\dagger(g)$ send the $n$-particle space into the
$(n+1)$--one; by this they are called annihilation and creation operators,
respectively. On $\Dscr$ the annihilation and creation operators satisfy the
\emph{canonical commutation rules} $\left[ a(g), a(f)\right]=0$, $\left[
a^\dagger(g), a^\dagger(f)\right]=0$, $\left[ a(g), a^\dagger(f)\right]= \langle
g|f\rangle_{{}_\Xscr}$, which characterize Bose fields (\cite{Par} pp.~136--146
and Example 24.1 p.~181). For the operators $A_j(t)$ and $A_j^\dagger(t)$ these
commutations relations become $ [A_j(t),A_i(s)] = 0$,
$[A_j^\dagger(t),A_i^\dagger(s)] = 0$, $[A_j(t),A_i^\dagger(s)] = \delta_{ij}
\min(t,s)$; in particular
\begin{equation}\label{d.8}
[A_j(t) + A_j^\dagger(t),A_i(s) + A_i^\dagger(s)] = 0\,, \qquad \forall i,j,t,s.
\end{equation}

A connection between creation and annihilation operators and Weyl operators is
given by the following theorem (\cite{Par} pp.~136, 139, 143--144).

\begin{theorem}\label{t2}
For every $k\in \Xscr$ the operator $a(k) + a^\dagger(k) $ is essentially
selfadjoint on $\Escr$. Moreover, if we denote again by $a(k) + a^\dagger(k) $
its selfadjoint extension, we have $\Wscr(\im k;\bbbone) = \exp\{\im [a(k) +
a^\dagger(k)]\}$.
\end{theorem}

>From now on, when $v\in \Zscr$, we shall set $v_j:= \langle e_j|v\rangle_\Zscr$.
Moreover, we shall confuse essentially selfadjoint operators with their
selfadjoint extensions.

Let us set $\Rscr := \left\{ k\in L^2_{\rm loc}(\mathbb R_+ ; \Zscr) : k_j(t)\in
\mathbb R, \, \forall j,t \right\}$. From eq.~\eqref{d.5} we have that
\begin{equation}\label{s.5}
\left[ \Wscr_{t_1} \left( \im k_1;\, \bbbone\right),\, \Wscr_{t_2}\left(\im k_2;
\, \bbbone\right)\right]=0\,, \qquad \forall t_1,\, t_2 \in \mathbb R_+\,,
\quad \forall k_1,\, k_2 \in \Rscr,
\end{equation}
or, which is the same, from the canonical commutation rules we have
\begin{equation}\label{s.6}
\left[ a\left( \chi_{[0,t_1]}k_1\right)+ a^\dagger \left( \chi_{[0,t_1]} k_1
\right),\,  a\left( \chi_{[0,t_2]}k_2\right)+ a^\dagger \left(
\chi_{[0,t_2]} k_2 \right)\right] = 0\,.
\end{equation}
Then, we set $\Mscr(t) := \big\{\Wscr_t(\im k;\bbbone),\, k\in \Rscr
\big\}^{\prime\prime}$ (double commutant in $\Lscr(\Gamma)$); $\Mscr(t)$ is the
von Neumann sub--algebra of $\Lscr(\Gamma)$ generated by the Weyl operators
$\big\{\Wscr_t(\im k;\bbbone),\, k\in \Rscr \big\}$ and it turns out to be a
commutative algebra. Moreover, by the spectral theorem, the previous one and the
definitions \eqref{d.6} and \eqref{d.7} we have also that every operator in
$\Mscr(t)$ can be seen as a bounded function of the commuting selfadjoint
operators $\big\{ A_j(s) + A_j^\dagger(s)\,,\ j=0,1,\ldots,\, s\in[0,t] \big\}$.

Let us consider now infinitely many independent standard Wiener processes
$W_j(t)$ as in Section 1, but let us choose a concrete probability space:
$\Omega$ is the set of all the functions $\omega : t \mapsto \omega(t) \equiv
\left( \omega_0(t), \omega_1(t), \ldots \right)$ with $\omega_j(0)=0$ and
$t\mapsto \omega_j(t)$ real and continuous, $\Fscr$ is the $\sigma$-algebra
generated by all canonical projections, $P$ is the Wiener measure on
$(\Omega,\Fscr)$, the Wiener processes are given by the canonical projections,
i.e.\ $W_j\big(t\big)(\omega) \equiv \omega_j(t)$, $\Fscr_t$ is the
sub-$\sigma$-algebra of $\Fscr$ generated by the family of random variables
$\big\{ W_j(s),\ j=0,1,\ldots,\ s\in [0,t]\big\}$. Let us call $L^2_{{}_W} :=
L^2(\Omega,\Fscr,P)$ the Wiener space; we recall that for $a,b \in L^2_{{}_W}$
one has
\begin{equation}\label{s.7}
\langle a|b\rangle_{L^2_W}= \int_\Omega \overline{ a(\omega)}\,
b(\omega)P(\dd \omega) = \mathbb E_{{}_P} \left[ \overline{a}\, b\right].
\end{equation}

\begin{theorem}\label{t3}
There exists a unique unitary isomorphism $\Iscr : \Gamma_0 \to L^2_{{}_W}$
satisfying
\begin{equation} \label{s.1}
{}\hskip-30pt\Iscr E(f) = E_{{}_W}(f):= \exp \biggl\{\sum_{j=0}^\infty
\int_0^{+\infty}\left[ f_j(t)\, \dd W_j(t) - \frac{1}{2} \left( f_j(t)
\right)^2 \dd t \right]\biggr\}, \qquad \forall f\in L^2(\mathbb R_+;\Zscr).
\end{equation}
Moreover, we have
\begin{equation}\label{s.3}
\Iscr\, \Wscr_t (\im k; \bbbone)\, \Iscr^{-1} = \exp \biggl\{ \im
\sum_{j=0}^\infty \int_0^t k_j(s)\,\dd W_j(s)\biggr\},\qquad
\forall k\in \Rscr, \ \forall t>0,
\end{equation}
and, $\forall j=0,1,\ldots$, $\forall t>0$,
\begin{equation}\label{s.4}
\Iscr \left( A_j(t) + A_j^\dagger(t) \right) \Iscr^{-1} = W_j(t)\,,
\qquad \Iscr\, \Mscr(t)\, \Iscr^{-1} = L^\infty(\Omega,\Fscr_t,P)\,;
\end{equation}
here $W_j(t)$ is a multiplication operator in $L^2_{{}_W}$.
\end{theorem}

\noindent {\it Proof.}
The isomorphism $\Iscr$ is constructed in \cite{Par} Example 19.9 pp.~130--131.

By eqs.~\eqref{d.4}, \eqref{d.2}, \eqref{s.1}, we have $\forall k\in
L^2_{\rm loc}(\mathbb R_+; \Zscr)$ and $\forall f\in L^2(\mathbb R_+; \Zscr)$
\begin{equation} \label{s.2}{}\hskip-20pt
\Iscr\, \Wscr_t(\im k;\bbbone)\, \Iscr^{-1}\, E_{{}_W}(f) = \exp \biggl\{
\im \sum_{j=0}^\infty \int_0^t \Bigl[ k_j(s) \,\dd W_j(s) +  \big( 2f_j(s)
+ \im k_j(s) \big) \Im k_j(s) \, \dd s \Bigr]\biggr\} E_{{}_W}(f).
\end{equation}
Then, Theorems \ref{t1}, \ref{t2} and eq.~\eqref{s.2} give eq.~\eqref{s.3} and
the first of eqs.~\eqref{s.4}. By eq.~\eqref{s.3} and the definition of
$\Mscr(t)$ we have also the second of eqs.~\eqref{s.4}.
\hfill{$\square$}

\smallskip

It is possible to extend the isomorphism $\Iscr$ to the whole $\Gamma \equiv
\Gamma_{-\infty}^0\otimes \Gamma_0$ by mapping it into $L^2_{{}_W}\otimes
L^2_{{}_W}$. Let us stress that $\Iscr$ is the unitary transformation which
simultaneously diagonalizes all the commuting selfadjoint operators $\{ A_j(t)
+A_j^\dagger(t)$, $j=0,1,\ldots$, $t\in \mathbb R_+\}$.

\section{Quantum stochastic differential equations}

Let $\Hscr$ be a separable complex Hilbert space as in Section 1. Hudson and
Parthasarathy \cite{HudP} developed a \emph{quantum stochastic calculus} which
gives meaning to integrals with respect to $\dd A_j(t)$, $\dd A_j^\dagger(t)$,
$\dd t$; we give only a rough idea of the definition, just what we shall need in
the following. Let $G(t)$ be an operator with domain including $\Hscr\otimes
\Escr$ $\big($the linear span of the vectors of the type $\varphi\otimes
E(f)\big)$ and such that $G(t) \, \Hscr\otimes \Escr_{-\infty}^t \subset \Hscr
\otimes \Gamma_{-\infty}^t$. Then, the integral $\int_0^t G(s)\dd A_j(s)$ is
defined as a suitable limit (if it exists) of a sum of terms of the type $G(t_i)
\left( A_j(t_{i+1}) - A_j(t_i) \right)$; this kind of definition is analogous to
the definition of classical It\^o's stochastic integral (\cite{Par} Sects.~24
and 25). For instance, one has
\begin{equation}\label{d.19}
a(k) = \sum_{j=0}^\infty \int_0^{+\infty} \overline{k_j(t)} \, \dd A_j(t)\,, \quad
a^\dagger(k) = \sum_{j=0}^\infty \int_0^{+\infty} k_j(t) \, \dd A_j^\dagger(t)\,,\quad
k\in L^2(\mathbb R_+;\Zscr).
\end{equation}

Let us introduce now the operators $H,\, L_j \in \Lscr(\Hscr)$, such that
$H=H^*$ and $\sum_{ j=0}^\infty L_j^* L_j$ strongly convergent in
$\Lscr(\Hscr)$; let us set $\widetilde K:= H - \frac{\im}{2} \sum_{j=0}^\infty
L_j^* L_j$. We consider the \emph{quantum stochastic Schr\"odinger equation}
\begin{align}\label{d.9}
\dd U(t,s) &= \biggl\{ \sum_{j=0}^\infty \left[ L_j \dd A_j^\dagger(t) - L_j^*
\dd A_j(t) \right] - \im \widetilde K\dd t \biggr\} U(t,s)\,, \\
U(s,s) &= \bbbone\,;
\label{d.10}\end{align}
here $s\leq t$, $U(t,s)$ is an operator on $\Hscr\otimes \Gamma$ and $L_j$ is
identified with $L_j \otimes \bbbone$, $A_j(t)$ with $\bbbone \otimes A_j(t)$
and so on.

\begin{theorem}\label{t4} {\rm (Hudson, Parthasarathy, Frigerio)}
Equations (\ref{d.9}), (\ref{d.10}) have a unique solution $U(t,s)$, $t\geq s$;
we have also $U(t,s) \in \Uscr(\Hscr \otimes \Gamma_s^t) \subset
\Uscr(\Hscr\otimes \Gamma)$ and $U(t,r)U(r,s)= U(t,s)$ for $s\leq r \leq t$. If
we fix $f\in \Xscr$ and set $\forall a \in \Lscr(\Hscr)$
\begin{equation}\label{d.20}
\widetilde \Lscr_t[a] := \im \bigg[ H+ \im \sum_{j=0}^\infty \left(
\overline{ f_j(t)}\, L_j - f_j(t) L^*_j \right), \,a \bigg] + \frac{1}{2}
\sum_{j=0}^\infty \left( \big[ L_j^*,a \big]L_j + L_j^* \big[ a, L_j \big] \right),
\end{equation}
then we have, $\forall \xi \in \Hscr$, $ \forall a \in \Lscr(\Hscr)$, $\forall
t,s \in \mathbb R$, $t \geq s$,
\begin{multline}\label{d.21}
\langle U(t,s) \xi \otimes e(f) | a \, U(t,s) \xi\otimes e(f)
\rangle_{{}_{\Hscr\otimes\Gamma}} = \langle \xi|a\xi\rangle_{{}_\Hscr} + {}
\\
{}+\int_s^t \langle U(\tau,s) \xi \otimes e(f) | \widetilde
\Lscr_\tau [a]\,  U(\tau,s) \xi\otimes e(f) \rangle_{{}_{\Hscr\otimes\Gamma}}\,\dd \tau.
\end{multline}
Moreover, the quantity $S^*(t+s) U(t+s,s) S(s)$ does not depend on $s$ and
$\{\widehat U(t),\ t\in \mathbb R\}$, where $\widehat U(t) :=  S^*(t) U(t,0)$
for $t\geq 0$ and $\widehat U(t):= U^*(-t,0) S(-t)$ for $t<0$, is a strongly
continuous one--parameter group of unitary operators on $\Hscr\otimes \Gamma$.
\end{theorem}

For the proof of the first part of the theorem see \cite{Par} Theorem 27.8
p.~228, Proposition 26.7 p.~216, Corollaries 27.9 and 27.10 p.~230; for the
statements about $\widehat U(t)$ see \cite{Fri1}. Let us stress that from the
strong continuity we have that $\widehat U(t)$ is strongly differentiable on
some dense domain and the same is true for $S(t)$. But the two infinitesimal
generators are unbounded and mutually related in such a way that $U(t,s) = S(t)
\widehat U(t-s) S^*(s)$ is not differentiable in the usual sense; however,
although the usual derivative of $U(t,s)$ does not exist, the quantum stochastic
calculus gives meaning to its stochastic differential \eqref{d.9}. Another
interesting point is that, from the connection between $\widehat U(t)$ and the
\emph{master equation} \eqref{d.21}, it is possible to deduce that the Stone
generator of $\widehat U(t)$ is a selfadjoint operator, necessarily unbounded
from above and from below.

In the standard formulation of quantum mechanics it is assumed that the
evolution of an isolated system is represented by a strongly continuous group of
unitary operators. So, by the previous theorem, we can think we have a quantum
system $A$, described in the Hilbert space $\Hscr$, interacting with a second
system $B$, a bosonic field described in $\Gamma$. The dynamics of $A+B$ is
given by $\widehat U(t)$; if we interpret $S(t)$ as the free dynamics of the
field, then $U(t,s)$ is the evolution operator of the compound system in the
\emph{interaction picture} with respect to the free--field dynamics. It is
possible to show \cite{Gar,QO} that evolutions such as $U(t,s)$ are a sensible
approximation to the dynamics of a photo--emissive source $A$, such as an atom, an optical
system in a cavity, $\ldots$, interacting with the electromagnetic field $B$.

Let us fix now $s=0$ as initial time and $\xi\otimes e(f)$ as initial state,
$\xi\in \Hscr$, $\|\xi\|_{{}_\Hscr} =1$, $f\in L^2_{\rm loc}(\mathbb R; \Zscr)$;
we have not defined coherent vectors for such an $f$, but we shall see in the
following that in our applications the queue of $f$ does not matter and that it
is often useful to consider periodic functions. Indeed the choice $f(t) \sim
\exp(-\im \omega_0 t)$ could represent a stimulating monochromatic laser.
Moreover, we set $U_t:= U(t,0)$, so that $U_t\, \xi \otimes e(f)$ is the state
of our compound system at time $t$.

Up to here we have studied the evolution of the system. Let us consider now a
measurement process in which some observables of the system are monitored with
continuity in time; the mathematical model of the detection scheme will be given
again in terms of operators on the Fock space and quantum stochastic integrals.
Quantum stochastic calculus has been introduced into the theory of continual
measurements in \cite{BLup}.

Let us recall that, in the usual formulation of quantum mechanics, observables
are represented by selfadjoint operators. Let us consider a quantum system
represented in a Hilbert space $\mathfrak h$ and let $\psi_t$ be the state of
the system at time $t$ ($\psi_t\in \mathfrak h$, $\|\psi_t\|_{\mathfrak h}=1$).
Let $X_1$, $\ldots$, $X_n$ be commuting selfadjoint operators with joint
spectral projections $\Pi_{\vec X}(B)$, $B$ Borel subset of $\mathbb R^n$. Then,
$\langle \psi_t| \Pi_{\vec X}(B)\psi_t\rangle_{\mathfrak h} $ is interpreted as
the probability that the measurement of the ``observable" $\vec X$ gives  a
result in $B$ at time $t$; indeed, $\langle \psi_t| \Pi_{\vec
X}(\cdot)\psi_t\rangle_{\mathfrak h} $ is a probability measure on $\mathbb
R^n$. Obviously, the whole probability measure is known if we know its Fourier
transform (the characteristic function), which turns out to be given by $\langle
\psi_t| \exp\{ \im  \vec k \cdot \vec X \}\, \psi_t\rangle_{\mathfrak h} $,
$\vec k \in \mathbb R^n$. We can also say  that we know completely the
probability law of $\vec X$ if we know all the quantities $\langle \psi_t|
a\,\psi_t\rangle_{\mathfrak h} $ with $a$ belonging to the commutative von
Neumann sub--algebra of $\Lscr(\mathfrak h)$ generated by the set of operators
$\{\exp[ \im  \vec k \cdot \vec X ]$, $\vec k \in \mathbb R^n\}$. This can be
generalized even to families of infinitely many commuting selfadjoint operators.

Now, let us consider again our system $A+B$ and the family of commuting
selfadjoint operators
\begin{equation}\label {o.1}
Z_j(t) := \sum_{i=0}^\infty \int_0^t \left[ \overline{ V_{ij}(s)}\, \dd A_i(s)
+ V_{ij}(s)\, \dd A_i^\dagger(s) \right], \qquad t\geq 0, \quad j=0,1,\ldots;
\end{equation}
here $V$, $V(t)$ are unitary operators as in eq.~(\ref{d.2}) and $ V_{ij}(t):=
\langle e_j|V(t) e_i\rangle_{{}_\Zscr}$. We consider $\{Z_j(t)\}$ as a set of
continuously monitored observables. Indeed, in the case of electromagnetic field emitted
by some source, the monitoring of observables of the type \eqref{o.1} can be concretely
realized by a measurement scheme called \emph{balanced heterodyne detection} \cite{QO}.
Let us denote by $\Mscr_Z(t)$ the von Neumann algebra generated by the Weyl operators
$\left\{ \exp \left[ \im
\sum_{j=0}^\infty \int_0^t k_j(s)\, \dd Z_j(s) \right], \ k\in \Rscr \right\}$.
By the previous discussion, the law of the measurement of $\vec Z(s)$ for $s$
from 0 to $t$ is known if we know the quantities $\langle U_t \xi \otimes e(f)|
Y U_t \xi \otimes e(f) \rangle_{{}_{\Hscr\otimes \Gamma_0}}$, $\forall Y\in
\Mscr_Z(t)$. The next step is to transform these quantities and to link them to
the stochastic processes introduced in Section 1.

>From now on $V$ is the unitary operator appearing in the observables
\eqref{o.1}, $f$ is the function characterizing the initial state of $B$, $\xi$
is the initial state of $A$ ($\xi\in \Hscr$, $\|\xi\|_{{}_\Hscr}=1$) and
$H_0=H_0^* \in \Lscr(\Hscr)$ is a selfadjoint operator which can be chosen
arbitrarily and is introduced for future convenience.

\begin{theorem}\label{t5} For every $t\in \mathbb R_+$ and $Y\in \Mscr_Z(t)$, we set
\begin{equation}\label{d.12}
\widetilde U_t := \exp[\im H_0 t]\, \Wscr_t^*(f;V) U_t \Wscr_t(f;V)\,,
\end{equation}
\begin{equation}\label{d.18}
\psi_t := \Iscr \,\widetilde U_t \,\xi\otimes e(0)\,,
\end{equation}
\begin{equation}\label{o.3}
\widetilde Y := \Wscr^*_t(f;V)\, Y\, \Wscr_t(f;V) \,.
\end{equation}
Then, we have $\widetilde U_0=\bbbone$, $\widetilde Y\in \Mscr(t)$, $\Iscr
\widetilde Y\Iscr^{-1}\in L^\infty(\Omega,\Fscr_t,P)$,
\begin{equation}\label{o.2}
\hskip -20pt \langle U_t \xi \otimes e(f)| Y U_t \xi \otimes e(f)
\rangle_{{}_{\Hscr\otimes \Gamma_0}}= \langle \widetilde U_t \xi
\otimes e(0)| \widetilde Y \widetilde U_t \xi \otimes e(0)
\rangle_{{}_{\Hscr\otimes \Gamma_0}} = \langle \psi_t |\Iscr
\widetilde Y \Iscr^{-1} \psi_t\rangle_{{}_{\Hscr\otimes L^2_W}}\,,
\end{equation}
\begin{equation}\label{d.13}
\dd \widetilde U_t =\biggl\{ \sum_{j=0}^\infty \left[  R_j(t)
\dd A_j^\dagger(t) - R_j^*(t) \dd A_j(t) \right] - \im K(t)\dd t \biggr\} \widetilde U_t\,,
\end{equation}
\begin{equation}\label{d.14}
K(t) := \e^{\im H_0 t} \bigg\{ \widetilde K - H_0 + \im
\sum_{j=0}^\infty \left[ \overline {f_j(t)}\, L_j - f_j(t)
L_j^*\right]\bigg\} \e^{-\im H_0 t},
\end{equation}
\begin{equation}\label{d.15}
R_j(t) := \sum_{i=0}^\infty \left( V^*(t)\right)_{ji}\,
\e^{\im H_0 t} L_i\,\e^{-\im H_0 t}\,;
\end{equation}
finally, $\psi_t$ satisfies the SDE \eqref{sd1} with initial condition
$\psi_0(\omega)=\xi$, $\forall \omega \in \Omega$. Moreover, if we assume
$\displaystyle\esssup_{t\leq T} \|f(t)\|_{{}_\Zscr} < +\infty$ for all $T\in
\mathbb R_+$, then the coefficients \eqref{d.14} and \eqref{d.15} satisfy
Assumptions \ref{A1} and \ref{A4}.
\end{theorem}

\noindent {\it Proof}. Note that from eqs.~(\ref{d.4})--(\ref{d.2}) we have
\begin{equation}\label{d.11}
\Wscr_t(f;V)e(0) = e\left(\chi_{[0,t]}f\right), \qquad
\Wscr\left(\left(1-\chi_{[0,t]}\right)f;\,\bbbone\right) \Wscr_t(f;V)e(0) = e(f)\,.
\end{equation}
The operator $U_t$ acts not trivially only on $\Hscr\otimes \Gamma_0^t$,
while $\Wscr \left(\left(1-\chi_{[0,t]}\right)f;\,\bbbone\right)$ only on
$\Gamma^0\otimes \Gamma_t$; therefore, these two operators commute. Then, by
eqs.~(\ref{d.11}) and (\ref{d.12}) one has
\begin{equation}\label{d.16}
U_t \,\xi\otimes e(f) = \Wscr \left( \left( 1- \chi_{[0,t]}\right)f;\,\bbbone\right)
\Wscr_t(f;V)e^{-\im H_0 t} \, \widetilde U_t \,\xi\otimes e(0)\,.
\end{equation}
The first of eqs.~\eqref{o.2} follows from eqs.~\eqref{d.16}, \eqref{o.3} and
the second one from eq.~\eqref{d.18}. The fact that $\widetilde Y$ belongs to
$\Mscr(t)$ follows from the definitions of $\Mscr(t)$ and $\Mscr_Z(t)$ and from
eqs.~\eqref{d.5} and \eqref{d.2}; the fact that one has $\Iscr \widetilde  Y
\Iscr^{-1} \in L^\infty(\Omega,\Fscr_t,P)$ follows from the second of
eqs.~\eqref{s.4}.

>From eqs.~\eqref{d.5}, \eqref{d.2}, \eqref{d.6}, \eqref{d.7} one obtains how the
operation $\Wscr_t^*(f;V)$ $\cdot$ $\Wscr_t(f;V)$ transforms the operators
$A_j(t)$, $A_j^\dagger(t)$; then, eqs.~\eqref{d.13}--\eqref{d.15} follow from
eqs.~\eqref{d.12}, \eqref{d.9}.

By eq.~(\ref{d.6}), the operators $A_j(t)$ annihilate the vacuum. Moreover, by
the factorization properties of $\Gamma$ and the definition of quantum
stochastic integrals, we have that $\dd A_j(t)$, acting on $\Gamma_t^{t+\dd t}$,
commutes with $\widetilde U_t$, acting on $\Hscr \otimes \Gamma_0^t$, or,
better, $\int_0^t \cdots \dd A_j(s) \widetilde U_s = \int_0^t \cdots  \widetilde
U_s \dd A_j(s)$. Therefore, when eq.~(\ref{d.13}) is applied to $\xi\otimes
e(0)$, the integral with respect to $A_j$ gives a vanishing contribution and the
integrand can be changed at will; so we can write \cite{Fri2}
\begin{equation}\label{d.17}
\dd \widetilde U_t \,\xi \otimes e(0) = \biggl\{ \sum_{j=0}^\infty R_j(t)\left[
\dd A_j^\dagger(t) + \dd A_j(t) \right] - \im K(t)\dd t \biggr\} \widetilde U_t
\,\xi \otimes e(0)\,.
\end{equation}
By the first of eqs.~\eqref{s.4} we have that $\psi_t$ satisfies the SDE
\eqref{sd1} with initial condition $\psi_0 = \xi\otimes E_{{}_W}(0)$; but
$E_{{}_W}(0) \equiv 1$, so that we can write $\psi_0(\omega)=\xi$, $\forall
\omega \in \Omega$.

Finally, it is possible to prove that $\left \| \sum_{j=0}^\infty R_j^*(t)
R_j(t) \right\| = \left \| \sum_{j=0}^\infty L_j^* L_j \right\|$ and $\|K(t)\|
\leq \|H- H_0 \| + \frac{1}{2} \left \| \sum_{j=0}^\infty L_j^* L_j \right\| +2
\|f(t)\| \, \left \| \sum_{j=0}^\infty L_j^* L_j \right\|^{1/2}$; then, one can
check that our coefficients satisfy Assumptions \ref{A1} and \ref{A4}.
\hfill{$\square$}

\smallskip

Let us comment the content of Theorem \ref{t5}. All the physical quantities
(probabilities, char\-ac\-ter\-is\-tic functional, moments) are given by the
``quantum expectations" $\langle U_t \xi \otimes e(f)| Y U_t \xi \otimes e(f)
\rangle_{{}_{\Hscr\otimes \Gamma_0}}$ with $Y\in \Mscr_Z(t)$. If we set
\begin{equation}\label{110}
\widehat Y := \Iscr \widetilde Y \Iscr^{-1} \equiv \Iscr \Wscr_t(f;V)Y \Wscr^*_t(f;V)
\Iscr^{-1}\,,
\end{equation}
we have that $\widehat Y$ is a random variable in $(\Omega, \Fscr_t, P)$ or
$\left(\Omega, \Fscr_t, \widehat P\right)$; moreover, by the definitions of $L^2_W$ and
$\widehat P$, we can write
\begin{equation}\label{111}
\left\langle \psi_t \left | \widehat Y \psi_t \right.\right\rangle_{{}_{\Hscr \otimes
L^2_W}} = \mathbb E_{{}_P} \left[ \left\| \psi_t \right\|^2_{{}_\Hscr} \widehat Y \right] =
\mathbb E_{{}_{\widehat P}} \left[  \widehat Y \right] .
\end{equation}
Finally, by eq.~\eqref{o.2} we obtain
\begin{equation}\label{112}
\langle U_t \xi \otimes e(f)| Y U_t \xi \otimes e(f) \rangle_{{}_{\Hscr\otimes \Gamma_0}} =
\mathbb E_{{}_{\widehat P}} \left[  \widehat Y \right],
\end{equation}
where $\widehat Y$ is given by eq.~\eqref{110}, $\widehat P$ by \eqref{sd2}, $\psi_t$
is the solution of \eqref{sd1} with non random initial condition $\xi$, $R_j(t)$
and $K(t)$ are given by eqs.~\eqref{d.14} and \eqref{d.15}. Therefore,
$\widehat P$ is the physical law of the output, as stated at the end of Section 1.

To discuss concrete physical applications of the previous formalism would be too
long. However, let us stress that the first non trivial example can be realized
in the Hilbert space $\Hscr=\mathbb C^2$. Indeed, in this space one can
construct a model describing a two--level atom of resonance frequency $\omega$
stimulated by a laser of frequency $\omega_0$ and emitting fluorescence light;
we are interested in the spectrum of the emitted light. Let us just give the
list of the choices that characterize the model: $\Hscr=\mathbb C^2$, $H= \omega
\begin{pmatrix} 1 & 0 \\ 0 & 0 \end{pmatrix}$ with $\omega > 0$, $L_j = \alpha_j
\begin{pmatrix} 0 & 0 \\ 1 & 0 \end{pmatrix}$ with $\alpha_j \in \mathbb C$ and
$0< \sum_j |\alpha_j|^2 < +\infty$, $f_j(t) = \lambda_j \, \e^{-\im \omega_0 t}$
with $\omega_0>0$, $\lambda_j \in \mathbb C$ and $\lambda_0=0$; the arbitrary
selfadjoint operator $H_0$ is chosen to be $H_0= \omega_0 \begin{pmatrix} 1 & 0
\\ 0 & 0 \end{pmatrix}$ (this eliminates some time dependence from the resulting
equations). The fluorescence light is made to interfere with a strong laser of
frequency $\nu$ and the intensity of the resulting light is measured with a
photodetector. It can be shown that this scheme corresponds to a continual
measurement of one of the components of $\vec Z$ (say $Z_0$) given by
eq.~\eqref{o.1} with $V_{ij}(s) = \e^{-\im \nu s}\, \delta_{ij}$, $\nu>0$. The
intensity of the emitted light can be obtained from the second moments $\mathbb
E_{{}_{\widehat P}} \left[ W_0(t_1) W_0(t_2) \right]$ of $W_0(t) = \Iscr
\Wscr_t(f;V) Z_0(t) \Wscr_t^*(f;V) \Iscr^{-1}$ ; as a function of $\nu$, this
intensity is the spectrum of the atom. Computations are given in \cite{Wigner}
where only quantum stochastic calculus is used (classical SDEs are not
explicitly introduced); the result is a three--peaked spectrum, already known in
the quantum--optical literature as the Mollow spectrum.

\end{document}